\documentclass{article}
\usepackage{amsbsy, amscd, amsfonts,amsgen, amsmath, amstext, amsxtra, amsthm, amsopn }
\newtheorem{theorem}{Theorem}

\newtheorem{lemma}{Lemma}
\newcommand{\al}{\alpha}
\newcommand{\be}{\beta}
\newcommand{\fr}{\frac}
\newcommand{\de}{\delta}
\newcommand{\ep}{\epsilon}
\newcommand{\e}{\epsilon}

\newcommand{\h}{\theta}
\def \th{\theta}
\def \w{\omega}
\def \C{C_{\be;02}}

\begin{document}

\title{Hausdorff Dimension and Hausdorff Measure for Non-integer based Cantor-type Sets}
\author{Qinghe Yin\\
{\it Mathematical Science Institute}\\
{\it The Australian National University, Australia}\\
Qinghe.Yin@maths.anu.edu.au}
\date{}
\maketitle
\begin{abstract}
We consider digits-deleted sets or Cantor-type sets with
$\be$-expansions. We calculate the Hausdorff dimension $d$ of these
sets and show that $d$ is continuous with respect to $\be$. The
$d$-dimentional Hausdorff measure of these sets is finite and
positive.
\end{abstract}

\section{Introduction}
The Hausdorff dimension and Hausdorff measure of expansions with
deleted digits are of interest to mathematicians.  The Cantor
middle-third set is a classical example. In 1993, M. Keane posed the
following question,see \cite{pollicott}:

\leftskip0.5truecm \noindent{\it Is the Hausdorff dimension
$d(\lambda)$ of the one parameter family of Cantor-type sets
$$
\Lambda(\lambda)=\left\{\sum_{k=1}^\infty i_k\lambda^k:\
i_k=0,1,3\right\}
$$
continuous for $\lambda\in[\frac14,\frac13]$?}

\leftskip0.truecm\noindent This question is answered by Pollicott
and Simon \cite{pollicott} in 1995. They show that
$d(\lambda)=\fr{\log3}{-\log\lambda}$ for almost all
$\lambda\in[\fr14,\fr13]$ and there exists a dense set with
$d(\lambda)<\fr{\log3}{-\log\lambda}$. Solomyak \cite{solomyak} in
1998 shows that the $d(\lambda)$-dimensional Hausdorff measure of
$\Lambda(\lambda)$ is zero for almost all $\lambda\in
[\fr14,\fr13]$. Keane, Smorodinsky and Solomyak study the size of
$\Lambda(\lambda)$ when $\lambda>\fr13$ (\cite{keaneetal}, 1995).

Notice the definition of $\Lambda(\lambda)$ implies no restriction
on the digits 0, 1 and 3. When $\lambda>\fr14$, different sequences
may express the same number. In this paper we study Hausdorff
dimension in a more restricted case, associated with
$\be$-expansions. We obtain very specific new results for
digits-deleted sets or Cantor-type sets with $\be$-expansions.

$\be$-transformations and $\be$-expansions are first introduced by
R\'enyi {\cite{Renyi}} in 1957 and further explored by Parry
{\cite{Parry}} in 1960. For fixed $\beta>1$, $T_{\beta}:\
[0,1)\rightarrow [0,1)$ is defined by
\begin{equation}
T_{\beta}x=\beta x-\lfloor\beta x\rfloor, \label{tbeta}
\end{equation}
where $\lfloor\cdot\rfloor$ is the floor function. By (\ref{tbeta})
we can define the $\beta$-expansion for any $x\in[0,1)$. Let
$a_1=\lfloor \beta x\rfloor$ and $a_n=\lfloor\beta
T_\beta^{n-1}x\rfloor$. Then
\begin{equation}
x=\frac{a_1}{\beta}+\frac{a_2}{\beta^2}+\cdots\label{ebeta}
\end{equation}
If we denote $x=(0.a_1a_2\cdots)_\beta$, then $T_\beta
x=(0.a_2\cdots)_\beta$. R\'enyi {\cite{Renyi}} and Parry
{\cite{Parry}} study the ergodic theory for $T_\be$.

Given $\be>2$, for
$0\le\h_0<\h_1<\cdots<\h_{q-1}\le\lfloor\be\rfloor$, we define
$$
C_{\be;\h_0\cdots\h_{q-1}}=
\overline{\left\{x=(0.a_1a_2\cdots)_\be:\
a_k\in\{\h_0,\h_1,\cdots,\h_{q-1}\}\right\}}.
$$
Corresponding to $\Lambda(\lambda)$ defined above, we have
$C_{\be;013}$ where $\be=\fr1\lambda\in(3,4]$. It is easy to see
that $\Lambda(\lambda)\supseteq C_{\be;013}$ with equality only when
$\lambda=\fr14$.   The classical Cantor middle-third set is
$C_{3;02}$. It is natural to ask the following questions:

\leftskip0.5truecm

 \noindent{\it What is is the Hausdorff dimension
and Hausdorff measure of $C_{\be;\h_0\cdots\h_{q-1}}$? Is
$\dim_H(C_{\be;\th_0\cdots\th_{q-1}})$ continuous with respect to
$\be$?}

 \leftskip0cm

Comparing $\Lambda(\lambda)$ and $C_{\be;013}$, we see that the
former is an attractor of the iterated function system (IFS) $(X;
f_0, f_1, f_2)$ 
where
$X=[0,\frac{3\lambda}{1-\lambda}]$ and
$$
f_0(x)=\lambda x,\ \ \ f_1(x)=\lambda(x+1)\ \text{ and }\
f_2(x)=\lambda(x+3),
$$
but when $\be$ is not an integer it is impossible  to consider
$C_{4;013}$ as an attractor of an IFS without introducing
elaborative additional machineries.

Given an IFS $(X; f_0,\cdots, f_{n-1})$ and $\be<n$, we  define the
$\be$-attractor as a certain compact subset of the attractor
determined by the $\be$-shift. We show that when all the maps are
similarities with the same contraction ratio $0<r<1$ and with a
separation condition, the Hausdorff dimension of the $\be$-attractor
is given by $s= \fr{\log\be}{-\log r}$ and the $s$-dimensional
Hausdorff measure is positive and finite.

For fixed $\be>2$ and the set of digits
$0\le\h_0<\cdots<\h_{q-1}\le\lfloor\be\rfloor$, we show that there
exists a number $\al>1$ such that $C_{\be,\h_0\cdots\h_{q-1}}$ is an
$\al$-attractor of an IFS.  The Hausdorff dimension of
$C_{\be,\h_0\cdots\h_{q-1}}$ is $s=\frac{\log\al}{\log\be}$ and the
$s$-dimensional Hausdorff measure of $C_{\be,\h_0\cdots\h_{q-1}}$ is
positive and finite. If the separation condition holds we obtain
this immediately from the general result for $\be$-attractors.
However the separation condition does not hold in general. We need a
direct proof. We also show that
 $\dim_H(C_{\be,\h_0\cdots\h_{q-1}})$ is
continuous with respect to $\be$ for $\be>\h_{q-1}$ and it has a
negative derivative with respect to $\be$ for almost all
$\be>\h_{q-1}$. There exists a nowhere dense subset of $(\h_{q-1},
\h_{q-1}+1]$ with Lebesgue measure 0 such that
$\dim_H(C_{\be,\h_0\cdots\h_{q-1}})$ has infinite derivative or
infinite one-sided derivative.

$\beta$-attractors fit into the more general {\it symbolic
construction} of Pesin and Weiss {\cite{pesin}}. A $\beta$-attractor
is the limit set of a symbolic construction using the $\beta$-shift.
If a symbolic construction is {\it regular}, Pesin and Weiss in
{\cite{pesin}} give a lower bound to the Hausdorff dimension of the
limit set using topological pressure. If $\be$ is simple then the
symbolic construction using the $\beta$-shift is regular under a
separation condition. It is not known this is true for all $\be>1$.
Barreira {\cite{bar}} gives a sufficient condition for a symbolic
construction is regular. However the construction of
$C_{\be,\h_0\cdots\h_{q-1}}$ does not obey this condition for some
value of $\beta$.
 If the symbolic construction is regular and the
 equilibrium measure associated with an $l$-estimating vector
  is a Gibbs measure then {\cite{pesin}}
 shows that the $s$-dimensional Hausdorff measure of the limit
set of a symbolic construction is positive, where $s$ is a lower
bound  of the Hausdorff dimension of the limit set related to the
$l$-estimating vector. It is known that if a symbolic system has the
specification property then any equilibrium measure is Gibbs
(\cite{ru}). A $\beta$-shift has the specification property if and
only if the length of strings of 0's in the $\beta$-expansion of 1
is bounded ({\cite{BM}}). Obviously this set does not contain any
interval. Schmeling {\cite{sch}} showed that this set has Hausdorff
dimension 1. Again, it is not known if every $\beta$-shift possesses
a Gibbs measure as its equilibrium measure associated with the
$l$-estimating vector of {\cite{pesin}}.

Although the symbolic construction for $\be$-attractors may not be
regular in general,  the symbolic construction of
$C_{\be,\h_0\cdots\h_{q-1}}$ is regular. It follows that the
Hausdorff dimension of $C_{\be,\h_0\cdots\h_{q-1}}$ can be obtained
from the results of \cite{pesin} by considering the topological
pressure of the corresponding $\al$-shift. However, as mentioned
above, in this way we can not get the information of Hausdorff
measure of $C_{\be,\h_0\cdots\h_{q-1}}$ in all cases. To overcome
this difficulty we give a new proof.

It is well known that the classical Cantor set can be defined to be
the digit-deleted set, $C=C_{3;02}$. It can also be defined
geometrically as the attractor of the IFS $F=\{[0,1];f_0,f_2\}$,
where $f_i(x)=\frac{x+i}3$ for $i=0,1,2$ are the inverse branches of
$T_3$. But for $2<\be<3$ the digit-deleted set $C_{\be;02}$ is not
the attractor of an IFS of inverse branches of $T_\be$ because of
the piecewise definition of $T^{-1}_\be$, namely
$$
T_\beta^{-1}x=\begin{cases}\frac x\beta,\ \ \frac {x+1}\beta,\ \
\frac{x+2}\beta,
\ \text{ for }x<\beta-2\\
\frac x\beta,\ \ \frac {x+1}\beta,\ \text{ for }x\ge\beta-2.
\end{cases}
$$
However we will show that there is an associated {\it local} IFS
(see {\cite{mbarnsley}}, p177) whose invariant sets are related to
$C_{\be;02}$ in an interesting way. A local IFS has its domain of at
least one of its functions not equal to the whole of the underlying
space, $[0,1]$ in the present case. We define
$$
F_\be=\{[0,1]; f_0: [0,1]\mapsto [0,1], f_2:[0,\be-2]\mapsto[0,1]\}
$$
where now $f_0(x)=\frac x\be$, and $f_2(x)=\frac{x+2}\be$. An
invariant set $A$ of the local IFS $F_\be$ is a nonempty compact
subset of $[0,1]$ such that
$$
A=f_0(A)\cup f_2(A\cap[0,\be-2]).
$$
In general, a local IFS may have no or many invariant sets
(\cite{mbarnsley}). We show that, with the exception of countably
many values of $\be$, $C_{\be;02}$ is the unique invariant set of
the local IFS $F_\be$. Otherwise, $F_\be$ possesses another
invariant set $B_{\be;02}$ which can be constructed by an interval
removal process. We have $C_{\be;02}\subset B_{\be;02}$ and
$B_{\be;02}\backslash C_{\be;02}$ consists of countably many
isolated points.

$\beta$-transformations and $\be$-expansions are of interest to
mathematicians in a broad range. After the pioneer works of R\'enyi
\cite{Renyi} and Parry {\cite{Parry}}, many research works related
to $\be$-transformations and $\be$-expansions have been published.
Among these works, \cite{Aki1}, \cite{Aki2}, \cite{Brown-Yin2},
\cite{Fan-Zhu}, \cite{Ito}, \cite{Ito2}, \cite{Luz} and \cite{Mori},
for example, have studied fractals or fractal sets related to
$\beta$-expansions. Barnsley \cite{barnsley} used
$\be$-transformation to study fractal tops.

 In section 2, we recall concepts and basic results
for Hausdorff dimension, iterated function systems, $\be$-expansions
and symbolic dynamical systems. In section 3, we study the Hausdorff
dimension and Hausdorff measure for $\be$-attractors under a
separation  condition. In section 4, we study the Hausdorff
dimension and Hausdorff measure for $C_{\be;\h_0\cdots\h_{q-1}}$,
the Cantor-type sets constructed by $\be$-expansions. In section 5,
we study the invariant sets of the local IFS $F_\be$.  In section 6,
we point out some interesting topics for further research.


\section{Hausdorff Dimension, Iterated Function Systems,
$\be$-expansion and Symbolic Dynamical Systems}

In this section we introduce concepts and definitions for Hausdorff
dimension, iterated function systems, $\be$-expansion and symbolic
dynamical systems.

Readers can find concepts of Hausdorff measure, Hausdorff
dimension and iterated function systems in \cite{Falconer} or
\cite{Falconer97}. For the original source of iterated function systems
see \cite{Hutchinson}.

Let $E$ be a subset of a metric space $X$. A $\delta$-cover of $E$
is a countable or finite collection $\{U_i\}$ of subsets of $X$ with
$|U_i|\le \delta$ such that $E\subset \cup_{i=1}^\infty U_i$, where
$|\cdot|$ is the diameter of the given set. For any $\delta>0$ and
$s\ge0$, define
$$
\mathcal{H}_\delta^s(E)=\inf\left\{\sum_{i=1}^\infty|U_i|^s:\
\{U_i\} \text{ is a $\delta$-cover of }E\right\}.
$$
Clearly, $\mathcal{H}_\delta^s(E)$ is decreasing with respect to
either $s$ or $\delta$.  Let
$$
\mathcal{H}^s(E)=\lim_{\delta\to0}\mathcal{H}_\de^s(E).
$$
We call $\mathcal{H}^s(E)$ the $s$-dimensional Hausdorff out measure
of $E$. The Hausdorff dimension of $E$ is defined by
$$
\dim_H(E)=\inf\{s\ge0:\ \mathcal{H}^s(E)=0\}=\sup\{s\ge0: \
\mathcal{H}^s(E)=\infty\}.
$$

An iterated function system (IFS) $F=(X; f_0, f_1, \cdots, f_{n-1})$
on a compact metric space $X$ consists of a number of contractions
$f_i: X\to X$, where $n\ge2$. There exists a non-empty compact
subset $E$ of $X$ such that
$$
E=\bigcup_{i=0}^{n-1}f_i(E).
$$
We say $E$ is the attractor of the IFS. For any sequence $i_1,
i_2,\cdots$ with $0\le i_k\le n-1$, and any $x\in X$, we have
$$
\lim_{k\to\infty}f_{i_1}\circ \cdots\circ f_{i_k}(x)\in E.
$$
The above limit exists and is independent of $x\in X$. On the other
hand, for any $y\in E$ there exists a sequence $(i_1,i_2,\cdots)$
such that
$$
y=\lim_{k\to\infty}f_{i_1}\circ \cdots\circ f_{i_k}(x)
$$
for any $x\in X$. We say $(i_1,i_2,\cdots)$ is an address of $y$.
One point in $E$ may have more than one addresses.

We say the Open Set Condition holds for an IFS, if there exists a
non-empty open set $O$ such that $f_i(O)\subset O$ and $f_i(O)\cap
f_j(O)=\emptyset$ for all $0\le i, j\le {n-1}$ and $i\ne j$.

Assume that $X$ is a compact subset of $\mathbf{R}^m$. If $f_i$ is a
similarity for all $i$ , i.e., for all $x, y\in X$,
$d(f_i(x),f_i(y))=r_id(x,y)$ for some $r_i<1$, and the Open Set
Condition holds, then the Hausdorff dimension of the attractor $E$
is given by
$$
\dim_H(E)=s
$$
where $s$ satisfies $r_0^s+r_1^s+\cdots+r_{n-1}^s=1$.

Let $\Sigma_n=\{0,1,\cdots,n-1\}^{\Bbb{N}}$. Then $\Sigma_n$ is
compact with respect to the product topology. Define a shift map
$\sigma:\Sigma_n\to\Sigma_n$ by
$$
\sigma(i_1,i_2,\cdots)=(i_2,i_3,\cdots).
$$
We call $\Sigma_n$ a full shift. A compact subset $\Sigma$ of
$\Sigma_n$ is said to be a subshift if $\sigma(\Sigma)=\Sigma$. A
block $(i_1,\cdots,i_k)$ is a forbidden word of $\Sigma$ if it does
not appear in any element of $\Sigma$. $\Sigma$ is a subshift of
finite type if it is determined by a finite set of forbidden words
(see \cite{lind}). Given an $n\times n$ 0-1 matrix $M$, let
$\Sigma_M$ contain all the elements $(i_1,i_2,\cdots)\in \Sigma_n$
such that $m_{i_k+1,i_{k+1}+1}=1$ for all $k\ge1$. We say $\Sigma_M$
is a Markov shift. Clearly, a Markov shift is a subshift of finite
type.

Fix $\be>1$. Assume that the $\be$-expansion of
$\be-\lfloor\be\rfloor$ is
$$
\be-\lfloor\be\rfloor=\fr{\ep_2}\be+\fr{\ep_3}{\be^2}+\cdots
$$
Let $\ep_1=\lfloor\be\rfloor$. Then
\begin{equation}
1=\fr{\ep_1}\be+\fr{\ep_2}{\be^2}+\fr{\ep_3}{\be^3}+\cdots\label{ex1}
\end{equation}
We say that (\ref{ex1}) is the $\be$-expansion of 1 and denote by
$1=(0.\e_1\e_2\cdots)_\be$. We say $\be$ is simple if the
$\be$-expansion of 1 has finite many non-zero terms.

Let $e_i=\e_i$ if $\be$ is non-simple;  or let $e_{kn+i}=\e_i$  and
$e_{(k+1)n}=\e_n-1$, for all $k\ge0$ and $i=1,2,\cdots, n-1$, if
$\be$ is simple with $\e_n>0$ and $\e_j=0$ when $j>n$. Given a
sequence $(a_1, a_2,\cdots)$ with
$a_i\in\{0,1,\cdots,\lfloor\be\rfloor\}$, the expression
\begin{equation}
\fr{a_1}{\be}+\fr{a_2}{\be^2}+\cdots\label{ex2}
\end{equation}
is the $\be$-expansion of some $x\in[0,1)$ if and only if for any
$n\ge1$ one has
$$
(a_{n},a_{n+1},\cdots)<(e_1,e_2,\cdots),
$$
where ``$<$" is according to the lexicographical order. When
(\ref{ex2}) is the $\be$-expansion for some $x\in[0,1)$, we say the
finite sequence $(a_1,a_2,\cdots,a_n)$ is $\be$-admissible.

Let
$$
1=\fr{\ep_1}\be+\fr{\ep_2}{\be^2}+\fr{\ep_3}{\be^3}+\cdots
$$
and
$$
1=\fr{\ep_1'}\al+\fr{\ep'_2}{\al^2}+\fr{\ep_3'}{\al^3}+\cdots
$$
be the $\be$ and $\al$-expansions of 1. If $\be<\al$ then
$$
(\ep_1,\ep_2,\cdots)<(\ep'_1,\ep_2',\cdots).
$$

Use $\Sigma_\be$ to denote the closure of all $\be$-admissible
sequences under product topology. Then $\Sigma_\be$ is a subshift
which we call the $\be$-shift. IF $\be=n$ is an integer, then
$\Sigma_\be$ is the full shift. If $\be$ is simple, then
$\Sigma_\be$ a subshift of finite type. For example, when
$\be=\fr{\sqrt5+1}2$, the golden
mean, $\Sigma_\be$ is a Markov shift with $M=\begin{pmatrix}1& 1\\ 1& 0\\
\end{pmatrix}$.


\section{$\be$-attractors and Their Hausdorff Dimension.}

Let  $(X; f_0,f_1,\cdots,f_{n-1})$ be an IFS with attractor $E$.
Define $\phi: \Sigma_n\mapsto E$ by
$$
\phi(i_1,i_2,\cdots)=\lim_{k\to\infty}f_{i_1}\circ\cdots \circ
f_{i_k}(x).
$$
$\phi$ is well defined since the limit in the right hand side
exists, is independent of $x$, and is in $E$.

\noindent{\bf Definition.} Given an IFS $(X;
f_0,f_1,\cdots,f_{n-1})$  with attractor $E$. For a subshift
$\Sigma$ of $\Sigma_n$ let
$$
E_\Sigma=\phi(\Sigma).
$$
We call $E_\Sigma$ the $\Sigma$-attractor of the given IFS. In
particular, when $\Sigma$ is a $\be$-shift for some $\be\le n$, we
use $E_\be$ to denote $E_{\Sigma_\be}$ and call it the
$\be$-attractor.

\noindent{\bf Remark.} The $\Sigma$-attractor defined above is a
special case of the more general symbolic construction of
\cite{pesin}. Many different settings of fractals with iterated
function systems can be viewed as $\Sigma$-attractors for some
subshifts. For example, the fractals defined by Markov shifts (see
\cite{bandt} or \cite{yinq}), the graph-directed fractals first
proposed by Mauldin and Williams \cite{mauldin} can fit into this
setting. However, when $\be$ is non-simple, it is difficult to fit
$E_\be$ into known settings other than the symbolic construction of
\cite{pesin}.

Recall that when all $f_i$'s are similarities with scale $r_i<1$ and
with the open set condition, the Hausdorff dimension of the
attractor $E$ is determined by $\sum_{i=0}^{n-1}r_i^s=1$. When all
the $r_i$'s are equal ($=r$, say), we have $s=\frac{\log n}{-\log
r}$.  For $\be$-attractors we have the following result.

\begin{theorem} Let $(X; f_0,f_1,\cdots,f_{n-1})$ be an IFS
such that
 $$
 d(f_i(x), f_i(y))=rd(x,y)
 $$
 for all $x, y\in X$ and $0\le i\le n-1$, where $0<r<1$.
Let $1<\be<n$. Assume that the $\be$-attractor $E_\be$ has the
following separation condition:
 $$
 f_i(E_\be)\cap f_j(E_\be)\cap E_\be=\emptyset, \ \text { for }i\ne j.
 $$
 Then the Hausdorff dimension of $E_\be$ is given by
$$
\dim(E_\be)=\frac{\log\be}{-\log r}.
$$
The  $s$-dimensional Hausdorff measure of $E_\be$ is positive and
finite, where $s=\dim_H(E_\be)$. \label{betaattractor}
\end{theorem}

Obviously, theorem {\ref{betaattractor}} holds when $\be$ is an
integer. Comparing the integer case, we may interpret the
$\be$-attractor as it is constructed by $\be$ many similarities. For
some value of $\be$, $\dim_H(E_\be)$ can be computed by existing
method. For example, if $\be=\frac{\sqrt5+1}2$, then $\Sigma_\be$ is
a Markov shift with $M=\begin{pmatrix}1&1\\1&0\end{pmatrix}$ (the
golden-mean shift). With the separation condition, the Hausdorff
dimension can be calculated by $\|MR^s\|=1$, where
$MR^s=\begin{pmatrix}r^s&r^s\\r^s&0\end{pmatrix}$ and $\|\cdot\|$ is
the Perron-Frobenus eigenvalue of the given matrix. Then we have
$1-r^s-r^{2s}=0$. This gives $r^s=\be^{-1}$ and $s=\fr{\log
\be}{-\log r}$. If $\be$ is given by $1=\fr1\be+\fr1{\be^3}$, then
the $\be$-expansion of real numbers induce a shift of finite type
with forbidden words determined by the set $\{(0,1,1),
(1,0,1),(1,1,0),(1,1,1)\}$. Define an IFS formed by the composed
maps, $\{f_0\circ f_0, f_0\circ f_1, f_1\circ f_0\}$. Then $E_\be$
is the Markov attractor with
$$
M=\begin{pmatrix} 1&1&1\\1&0&0\\1&1&0\end{pmatrix}.
$$
Noting that the  composed maps are similarities with scales $r^2$,
then with the separation condition, the Hausdorff dimension of
$E_\be$ is given by $\|MR^{2s}\|=1$. This gives us
$$
r^{6s}+2r^{4s}+r^{2s}-1=0.
$$
But
$$
\lambda^6+2\lambda^4+\lambda^2-1=(\lambda^3+\lambda-1)(\lambda^3+\lambda+1).
$$
Hence we have $r^{3s}+r^s=1$ which implies $r^s=\be^{-1}$, and
therefore $s=\frac{\log \be}{-\log r}$.

Obviously, this discussion can only apply to particular values of
$\be$ when the set of forbidden words are ``short" and the
transition matrix is ``small" of size whose eigenvalue is
``calculatable". It is hard to apply to general case. Besides, it is
not applicable when the related shift is not a subshift of finite
type. 

Use $\mathcal{S}_\be^k$ to denote the set of all $\be$-admissible
sequences of length $k$.  To prove Theorem \ref{betaattractor} we
need to estimate the size of $\mathcal{S}_\be^k$. The following
result can be found in \cite{Renyi} (equations (4.9), (4.10)).

\begin{lemma}  We have
$$
\be^k\le|\mathcal{S}_\be^k|\le\fr{\be^{k+1}}{\be-1}.
$$
 \label{lemma1}
\end{lemma}

\noindent{\it Proof of Theorem \ref{betaattractor}}. First we show
that $\dim_H(E_\be)\le \fr{\log\be}{-\log r}$.  For $(i_1,\cdots,
i_k)\in\mathcal{S}_\be^k\}$, denote $\Delta_{i_1\cdots
i_k}=f_{i_1}\circ \cdots\circ f_{i_k}(X)$. Then the collection
$\{\Delta_{i_1\cdots i_k}: (i_1,\cdots, i_k)\in\mathcal{S}_\be^k\}$
is a cover of $E_\be$. Since all $f_i$'s are similarities with scale
$r$, we have $|\Delta_{i_1\cdots i_k}|\le r^k|X|$. By Lemma
\ref{lemma1}, we get
$$
\sum_{{i_1\cdots i_k}\in\mathcal{S}_\be^k} |\Delta_{i_1\cdots
i_k}|^s \le
\fr{\be^{k+1}}{\be-1}r^{sk}|E_\be|^s=\fr{\be}{\be-1}|X|^s,
$$
where $s=\frac{\log\be}{-\log r}$. This shows that $\dim(E_\be)\le
s$ and the $s$-dimensional Hausdorff measure of $E_\be$ is finite.

Next we show that $\dim(E_\be)\ge \fr{\log\be}{-\log r}$. We will
show that for some $\delta_0>0$ we have
$\mathcal{H}_{\delta_0}^s(E_\be)>0$ where $s=\frac{\log\be}{-\log
r}$. Then $\mathcal{H}^s(E_\be)=\lim_{\delta\to0}
\mathcal{H}_\delta^s(E_\be)\ge\mathcal{H}_{\delta_0}^s(E_\be)$. This
gives $\dim_H(E_\be)\ge s$ and the $s$-dimensional Hausdorff measure
of $E_\be$ is positive.

 By the separation condition,
$$
\delta_0=\min_{i\ne j}\{d(f_i(E_\be)\cap E_\be,f_j(E_\be)\cap
E_\be\}>0.
$$
Then there exists $l>0$ such that $r^{l+1}|X|<\delta_0\le r^l|X|$.
Let $\mathcal{U}=\{U_i\}$ be a $\delta_0$ cover of $E_\be$. Since
$E_\be$ is compact without loss of generality we may assume that
$\mathcal{U}=\{U_1,\cdots,U_N\}$ is a finite cover. We may also
assume that $U_i\subset E_\be$. In fact we can use $U_i\cap E_\be$
to replace $U_i$. Choose $k_i$ such that
$$
r^{k_i+1}|X|<|U_i|\le r^{k_i}|X|.
$$
We claim that $U_i\subset\Delta_{i_1\cdots i_{k_i-l}}$ for some
$\be$-admissible $(i_1,\cdots, i_{k_i-l})$. Taking $y\in U_i$ which
has an address $(i_1,i_2,\cdots)\in\Sigma_\be$, for any $z\in E_\be$
with an address $(j_1, j_2,\cdots)\in\Sigma_\be$ such that
$(i_1,\cdots, i_{k_i-l})\ne(j_1,\cdots, j_{k_i-l})$. Assume that
$i_1=j_1,\cdots, i_p=j_p$ but $i_{p+1}\ne j_{p+1}$ for some
$p<k_i-l$. Then
$$
d(y, z)=d(f_{i_1}\circ\cdots\circ f_{i_{p}}\circ f_{i_{p+1}}(y_1),
f_{i_1}\circ\cdots\circ f_{i_{p}}\circ f_{j_{p+1}}(z_1))
$$
for some $y_1, z_1\in E_\be$. It is obvious that $f_{i_{p+1}}(y_1),
f_{j_{p+1}}(z_1)\in E_\be$. Then,
$$
\begin{aligned}
&d(y,z)=r^{p}d(f_{i_{p+1}}(y_1), f_{j_{p+1}}(z_1))\\ \ge
&r^{p}d(f_{i_{p+l}}(E_\be)\cap E_\be, f_{j_{p+l}}(E_\be)\cap E_\be) \ge r^{p}\delta_0\\
>&r^{p}r^{l+1}|X|\ge r^{k_i-l-1}r^{l+1}|X|\\
=&r^{k_i}|X|>|U_i|.
\end{aligned}
$$
Hence $z\notin U_i$ and $U_i\subset\Delta_{i_1\cdots i_{k_i-l}}$.
Now we get another cover of $E_\be$:
$$
\mathcal{C}=\{\Delta_{i_1\cdots i_{k_i-l}}:\ U_i\subset
\Delta_{i_1\cdots i_{k_i-l}},\ r^{k_i+1}|X|<|U_i|\le r^{k_i}|X|,\
U_i\in \mathcal{U}\}.
$$
For  $s>0$ we have
\begin{equation}
\sum_{\Delta_{i_1\cdots i_{k_i-l}}\in\mathcal{C}} |\Delta_{i_1\cdots
i_{k_i-l}}|^s=\sum_{\Delta_{i_1\cdots i_{k_i-l}}\in\Delta}
r^{s(k_i-l)}|X|^s\le r^{s(-l-1)}\sum_{U_i\in \mathcal{U}}|U_i|^s
\label{eqww}
\end{equation}
Let
 $$
 k= \max_{U_i\in \mathcal{U}}\{k_i-l:\ r^{k_i+1}|X|<|U_i|\le
r^{k_i}|X|\}.
$$
 Refine $\mathcal{C}$ into
$$
\mathcal{C}'=\{\Delta_{i_1\cdots i_k}:\ (i_1,\cdots,
i_k)\in\mathcal{S}_\be^k\}.
$$
Then $ |\Delta_{i_1\cdots i_{k}}|^s=r^{s(k-k_i)}|\Delta_{i_1\cdots
i_{ki-l}}|^s. $ Set $s=\fr{\log\be}{-\log r}$, then
$$
|\Delta_{i_1\cdots i_{k}}|^s=\be^{-(k-k_i)}|\Delta_{i_1\cdots
i_{k_i-l}}|^s.
$$
By Lemma \ref{lemma1} we see that each $\Delta_{i_1\cdots
i_{k_i-l}}\in \mathcal{C}$ contains at most
$\fr{\be^{k-k_i+l+1}}{\be-1}$ many $\Delta_{i_1\cdots i_k}\in
\mathcal{C}'$. Then
$$
\underset{j_1=i_1,\cdots,
j_{k_i-l}=i_{k_i-l}}{\sum_{(j_1,\cdots,j_k)\in\mathcal{S}_\be^k}}|\Delta_{j_1\cdots
j_k}|^s \le\fr{\be^{l+1}}{\be-1}|\Delta_{i_1\cdots i_{k_i-l}}|^s.
$$
Therefore,
$$
\sum_{\Delta_{i_1\cdots i_{k_i-l}}\in\mathcal{C}} |\Delta_{i_1\cdots
i_{k_i-l}}|^s \ge \be^{-l-1}(\be-1) \sum_{\Delta_{i_1\cdots
i_{k}}\in\mathcal{C}'}|\Delta_{i_1\cdots i_k}|^s
$$
The right hand side contains all the $\be$-admissible sequences in
$\mathcal{S}_\be^k$ and we have $|\Delta_{i_1\cdots i_k}|=r^k|X|$.
Hence
\begin{equation}
\begin{aligned}
&\sum_{\Delta_{i_1\cdots i_{k_i-l}}\in\mathcal{C}}
|\Delta_{i_1\cdots i_{k_i-l}}|^s
\ge\be^{-l-1}(\be-1)\sum_{(i_1,\cdots,
i_k)\in\mathcal{S}_\be^k}r^{sk}|X|^s \\
\ge&\be^{-l-1}(\be-1)
\be^kr^{sk}|X|^s=\be^{-l-1}(\be-1)|X|^s.\end{aligned} \label{eqwwq}
\end{equation}
 By (\ref{eqww}) and (\ref{eqwwq}) we obtain
$$
\begin{aligned}
&\sum_{U_i\in \mathcal{U}}|U_i|^s\ge
r^{s(l+1)}\sum_{\Delta_{i_1\cdots i_{k_i-l}}\in\mathcal{C}}
|\Delta_{i_1\cdots i_{k_i-l}}|^s\\
 \ge &r^{s(l+1)}\be^{-l-1}(\be-1)|X|^s
=\be^{-2(l+1)}(\be-1)|X|^s\end{aligned}
$$
for any finite $\delta_0$-cover $\mathcal{U}$. Hence
$\dim_H(E_\be)\ge s$ and $\mathcal{H}^s(E_\be)>0$. $\Box$


\section{Hausdorff Dimension for $C_{\be;\th_0\cdots\th_{q-1}}$}
\def \CC{C_{\be;\th_0\cdots\th_{q-1}}}
\def \w{\omega}
\def \aa{\al_{\be;\h_0\cdots\h_{q-1}}}

Let $C_{\be;\th_0\cdots\th_{q-1}}$ be as defined in section 1. Use
$z_\be$ to denote the maximum of $\CC$. We have the following
result.

\begin{lemma} $z_\be$ can be expressed as
\begin{equation}
z_\be=\fr{z_1}{\be}+\fr{z_2}{\be^2}+\cdots\label{zz}
\end{equation}
where $z_i\in\{\h_0,\cdots,\h_{q-1}\}$, and
\begin{equation}
(z_i,z_{i+1},\cdots)\le(z_1,z_2,\cdots).\label{zzz}
\end{equation}
\label{lz}
\end{lemma}

\noindent{\it Proof.} By the definition of $\CC$, we have two
possibilities: the digits of the $\be$-expansion of $z_\be$ are in
$\{\h_0,\cdots,\h_{q-1}\}$, or $z_\be$ is a limit of numbers whose
$\be$-expansion consists of digits in $\{\h_0,\cdots,\h_{q-1}\}$.
For the first case, the $\be$-expansion of $z_\be$ satisfies the
requirement of Lemma \ref{lz}. For the second case, there exits a
sequence $\{y_k\}$ with digits in $\{\h_0,\cdots,\h_{q-1}\}$ such
that $y_k\uparrow z_\be$. Assume that
$y_n=(0.y_1^{k}y_2^k\cdots)_\be$. Then we have
$$
(y_1^k, y_2^k,\cdots)<(y_1^{k+1}, y_2^{k+1}, \cdots).
$$
Since $\Sigma_\be$ is compact under the product topology, the
sequence $\{(y_1^k, y_2^k,\cdots)\}$ has a limit $(z_1,z_2,\cdots)$.
It is obvious that $(z_1,z_2,\cdots)$ satisfies (\ref{zz}). We show
that it  satisfies (\ref{zzz}) as well. If it is not true, then
$(z_k,z_{k+1},\cdots)>(z_1,z_2,\cdots)$. Let
$\tilde{y}_n=(0.y_k^ny_{k+1}^n\cdots)_\be$. Then
$$
\lim_{n\to\infty}\tilde{y}_n=\fr{z_k}{\be}+\fr{z_{k+1}}{\be^2}+\cdots>z_\be,
$$
a contradiction. \hskip10pt $\Box$

Using $(z_1, z_2, \cdots)$ we define a new sequence
$(\w_1,\w_2,\cdots)$ with $\w_i=j$ if $z_i=\h_j$. Then
$$
(\w_k,\w_{k+1},\cdots)\le(\w_1,\w_2,\cdots).
$$
Let $\aa$ be determined by
\begin{equation}
1=\fr{\w_1}{\al}+\fr{\w_2}{\al^2}+\cdots\label{alp}
\end{equation}
Then either (\ref{alp}) is the $\aa$-expansion of 1 or
$(\w_1,\w_2,\cdots)$ is periodic.

\begin{theorem} The Hausdorff dimension of $\CC$ is given by
$$
\dim_H(C_{\be;\th_0\cdots\th_{q-1}})=\fr{\log\al_{\be;\th_0\cdots\th_{q-1}}}{\log\be}.
$$
The $s$-dimensional Hausdorff measure of $\CC$  is positive and
finite where $s=\dim_H(\CC)$.\label{newtheo}
\end{theorem}

\noindent{\it Proof.} For $0\le i\le q-1$ define $f_i:
\left[0,\frac\be{\be-1}\right]\mapsto\left[0,\frac\be{\be-1}\right]$
by $f_i(x)=\frac {x+\th_i}\be$. We will show that
$C_{\be;\th_0\cdots\th_{q-1}}$ is the $\al$-attractor of the IFS
$\left(\left[0,\frac\be{\be-1}\right]; f_0,\cdots, f_{q-1}\right)$.
Here and rest of the proof we use $\al$ to denote
$\al_{\be;\th_0\cdots\th_{q-1}}$. Then by the proof of Theorem
{\ref{betaattractor}} we get that
$$
\dim_H(C_{\be;\th_0\cdots\th_{q-1}})\le\fr{\log\al}{\log\be}
$$
and the $s$-dimensional Hausdorff measure of $\CC$ is finite for
$s=\fr{\log\al}{\log\be}$. If we further have the separation
condition, then Theorem \ref{newtheo} is proved.

First we show that $\CC$ is an $\al$-attractor. Use
$\Sigma_{{\be;\th_0\cdots\th_{q-1}}}$ to denote the closure of all
sequences $(a_1,a_2,\cdots)$ with $a_i\in\{\h_0,\cdots,\h_{q-1}\}$
which forms the $\be$-expansion of some $x\in[0,1)$. Then
$$
(a_1,a_2,\cdots)\mapsto \fr{a_1}{\be}+\fr{a_2}{\be^2}+\cdots
$$
is a 1-1 correspondence between
$\Sigma_{{\be;\th_0\cdots\th_{q-1}}}$ and $\CC$.

Given $x\in C_{\be;\th_0\cdots\th_{q-1}}$, assume that
$x=\fr{a_1}{\be}+\fr{a_2}{\be^2}+\cdots$ for some
$(a_1,a_2,\cdots)\in\Sigma_{\be;\th_0\cdots\th_{q-1}}$. Define
$(i_1,i_2,\cdots)\in \Sigma_{\al}$ according to $a_k=\th_{i_k}$.
Then $x=\lim_{k\to\infty}x_k$ where
$$
x_k=\fr{a_1}{\be}+\fr{a_2}{\be^2}+\cdots+\fr{a_k}{\be^k}
=f_{i_1}\circ f_{i_2}\circ\cdots\circ f_{i_k}(0).
$$
On the other hand, for any $(i_1,i_2,\cdots)\in \Sigma_{\al}$ and
any $k>0$ we have\break
$(\th_{i_1},\th_{i_2},\cdots,\th_{i_k},\th_{i_0},\cdots)
\in\Sigma_{\be;\th_0\cdots\th_{m-1}}$. Denote
$x^*=\fr{\th_0}{\be}+\fr{\th_0}{\be^2}+\cdots$. Then
$$
\begin{aligned}
& x_k=f_{i_1}\circ f_{i_2}\circ\cdots\circ f_{i_k}(x^*)\\
=&\fr{\th_{i_1}}{\be}+\fr{\th_{i_2}}{\be^2}+\cdots
+\fr{\th_{i_k}}{\be^k}+\fr{\th_0}{\be^{k+1}}+\cdots \in
C_{\be;\th_0\cdots\th_{q-1}}.
\end{aligned}
$$
Hence $\lim_{k\to\infty}f_{i_1}\circ f_{i_2}\circ\cdots\circ
f_{i_k}(x^*)\in C_{\be;\th_0\cdots\th_{q-1}}$. This shows that
$C_{\be;\th_0\cdots\th_{q-1}}$ is the $\al$-attractor of the IFS
$\left(\left[0,\frac\be{\be-1}\right]; f_0,\cdots, f_{q-1}\right)$.

If we have $f_i(\CC)\cap f_j(\CC)=\emptyset$ for $i\ne j$ then
Theorem \ref{newtheo} can be obtained by Theorem
\ref{betaattractor}. Unfortunately in some cases we do not have the
separation condition. For example, when $\be=4$, the separation
condition does not hold for $C_{4;013}$. For Hausdorff dimension we
may use the following argument.

For $\al'<\al$, we use $\CC^{\al'}$ to denote the $\al'$-attractor
of the IFS in consideration. Then $\CC^{\al'}\subset \CC$. Since
$\al'<\al$, we have $1\notin\CC^{\al'}$. Then the separation
condition holds for $\CC^{\al'}$. Therefore,
$$
\dim_H(\CC)\ge\dim_H(\CC^{\al'})=\fr{\log\al'}{\log\be}
$$
for any $\al'<\al$ This implies that
$$
\dim_H(\CC)\ge\fr{\log\al}{\log\be}.
$$

However, this does not supply any information about the Hausdorff
measure of dimension $\fr{\log\al}{\log\be}$. We will use similar
discussion as the proof of Theorem \ref{betaattractor}.

Notice that when separation condition fails we must have $0, 1\in
\CC$. For an $\al$-admissible sequence $(i_1,\cdots,i_k)$ let
$I_{i_1\cdots i_k}$ to denote the smallest closed interval which
contains all numbers whose $\be$-expansion starting with
$(\th_{i_1},\cdots,\th_{i_k})$. If $(i_1,\cdots,i_k)\ne(j_1,\cdots,
j_k)$ then $I_{i_1\cdots i_k}\cap I_{j_1\cdots j_k}$ contains at
most one point. Let $\mathcal{U}=\{U_1,\cdots, U_n\}$ be a finite
cover of $\CC$. If $\be^{-k-1}<|U_i|\le\be^{-k}$ then $U_i$
intersects at most three $I_{i_1\cdots i_k}$. Then
$$
\sum_{I_{i_1\cdots i_k}\cap U_i\ne\emptyset}|I_{i_1\cdots i_k}|^s
\le 3\be^{-ks}<3\be^{s}|U_i|^s
$$
Assume that $|U_i|>\be^{-l}$ for all $i$. By Lemma 1
$(i_1,\cdots,i_k)$ can be extend to at most
$\frac{\al^{l-k+1}}{\al-1}$ many $(i_1,\cdots,i_k,\cdots,i_l)$. Then
$$
\sum_{I_{i_1\cdots i_l}\cap U_i\ne\emptyset}|I_{i_1\cdots i_l}|^s
\le\frac{3\al^{l-k+1}}{\al-1} \be^{-ls}=\frac{3\al^{-k+1}}{\al-1}
<\frac{3\al^{2}|U_i|^s}{\al-1}
$$
where $s=\frac{\log\al}{\log\be}$. Hence
\begin{equation}
\sum_{i=1}^N|U_i|^s>\frac{\al-1}{3\al}\sum_{(i_1,\cdots,i_l)}|I_{i_1\cdots
i_l}|^s\label{aha}
\end{equation}
where the sum is over all $\al$-admissible $(i_1,\cdots,i_l)$.

We claim that if the separation condition does not hold then we have
$\al>2$. In fact if $\al\le2$ we have the following three
possibilities:

1. $\th_1<\lfloor\be\rfloor$, which implies $1\notin\CC$;

2. $\th_0>0$, which implies $0\notin\CC$;

3. $\th_1\ge\th_0+1$.

\noindent In all these three cases we have $f_0(\CC)\cap
f_1(CC)=\emptyset$.

Now we estimate the number of $I_{i_1\cdots i_l}$ that
$|I_{i_1\cdots i_l}=\be^{-l}$. Let
$$
A_k=\{(i_1,\cdots,i_k)|(i_1,\cdots, i_{k-1},i_k+1) \text{ is
$\al$-admissible}\}.
$$
If $(i_1,\cdots, i_k)\in A_k$ then $(\th_{i_1},\cdots,
\th_{i_{k-1}},\th_{i_k+1})$ is $\be$-admissible. Thus\break
$(\th_{i_1},\cdots, \th_{i_{k-1}},\th_{i_k}+1)$ is $\be$-admissible
since $\th_{i_k}+1\le\th_{i_k+1}$. Therefore we have $|I_{i_1\cdots
i_k}|=\be^{-k}$. Let $\mathcal{S}_\al^k$ of all $\al$-admissible
$(i_1,\cdots, i_k)$. By Lemma 1 we have
$|\mathcal{S}_\al^k|\ge\be^k$. We show that $|A_k|\ge
c|\mathcal{S}_\al^k|$ for some constant $c$.

 Let
$B_k=\mathcal{S}_\al^k\backslash A_k$. We estimate
$\frac{|B_k|}{|A_k|}$. Notice that $|B_{k+1}|\le|\mathcal{S}_\al^k|
=|A_k|+|B_k|$ and $|A_{k+1}|\ge\lfloor\al\rfloor\cdot|A_k|$. Then
$$
\frac{|B_{k+1}|}{|A_{k+1}|}\le\frac{|A_k|+|B_k|}{\lfloor\al\rfloor\cdot|A_k|}
=\frac1{\lfloor\al\rfloor}\left(1+\frac{|B_k|}{|A_k|}\right).
$$
Continue this discussion we get that for any $k$ we have
$$
\frac{|B_k|}{|A_k|}\le\frac1{\lfloor\al\rfloor}+\frac1{\lfloor\al\rfloor^2}+\cdots
=\frac{1}{\lfloor\al\rfloor-1}.
$$
Therefore,
$$
\frac{|A_l|}{|\mathcal{S}_\al^l|}\ge
\frac1{1+\frac1{\lfloor\al\rfloor-1}}=\frac{\lfloor\al\rfloor-1}{\lfloor\al\rfloor}.
$$

By this and (\ref{aha}) we obtain
$$
\begin{aligned}
&\sum_{i=1}^N|U_i|^s>\frac{\al-1}{3\al}\sum_{(i_1,\cdots,i_l)\in
A_l}|I_{i_1\cdots i_l}|^s\\
\ge&\frac{\al-1}{3\al}\cdot\frac{\lfloor\al\rfloor-1}{\lfloor\al\rfloor}
\cdot|\mathcal{S}_\al^l|\cdot\be^{-sl}\\
\ge&\frac{\al-1}{3\al}\cdot\frac{\lfloor\al\rfloor-1}{\lfloor\al\rfloor}
 \cdot\al^l\be^{-sl}=\frac{\al-1}{3\al}\cdot\frac{\lfloor\al\rfloor-1}{\lfloor\al\rfloor}
\end{aligned}
$$
for all finite cover $\mathcal{U}$, where
$s=\frac{\log\al}{\log\be}$. This proves that the $s$-dimensional
Hausdorff measure of $\CC$ is positive. \hskip10pt $\Box$

{\noindent \bf Remark.} Observe that the Hausdorff dimension of
$C_{\be;\h_0\cdot\h_q}$ is related to ``base change". For the
classical Cantor set $C$, if we change base 3 into base 2 and the
digit 2 into 1, we get an (almost) 1-1 map from $C$ to $[0,1]$.
Theorem {\ref{newtheo}} demonstrates that the ``base change"
property also holds for non-integer bases.


\noindent{\bf Example.} Let $\be_1$, $\be_2$ be given by $
1=\fr{3}{\be_1}+\fr{2}{\be_1^2} $ and $
1=\fr{3}{\be_2}+\fr{3}{\be_2^2} $. Then $\be_1=\fr{3+\sqrt{17}}2$
and $\be_2=\fr{3+\sqrt{21}}2$. We consider $C_{\be;013}$ for
$\be\in[\be_1,\be_2]$.  It is easy to see that
$\Sigma_{\be_1;013}=\Sigma_{\be_2;013}$, and hence
$\Sigma_{\be;013}=\Sigma_{\be_1;013}$ for any $\be\in[\be_1,\be_2]$.
The maximal sequence of $\Sigma_{\be_1;013}$ is $(3,1,3,1,\cdots)$.
Thus $\al_{\be;013}$ is determined by
$$
1=\fr2\al+\fr1{\al^2}+\fr2{\al^3}+\fr1{\al^4}+\cdots
$$
which gives $\al=1+\sqrt3$. therefore, for any $\be\in[\be_1,\be_2]$
$$
\dim_H(C_{\be;013})=\fr{\log(1+\sqrt3)}{\log\be}.
$$
By this we see that $\dim_H(C_{\be;013})$ is decreasing with respect
to $\be$ on the interval $[\be_1,\be_2]$.

\begin{theorem} For fixed $\th_0, \cdots, \th_{q-1}$,
$\dim_H(C_{\be;\th_0\cdots\th_{q-1}})$ is continuous for
$\be>\th_{q-1}$. \label{cont}
\end{theorem}

\noindent{\it Proof.} We need to prove that
$\al_{\be;\th_0\cdots\th_{q-1}}$ is continuous with respect to $\be$
for $\be>\th_{q-1}$. Obviously, $\al_{\be;\th_0\cdots\th_{q-1}}$ is
non-decreasing with respect to $\be$. It is clear that we have
$\al_{\be;\th_0\cdots\th_{q-1}}>q-1$. If we can show that for any
$\al\in(q-1,q]$, there exists a $\be$ such that
$\al=\al_{\be;\th_0\cdots\th_{q-1}}$ by the monotone property we see
that it is continuous.

Given $\al\in(q-1,q]$, assume that the $\al$-expansion of 1 is
$$
1=\fr{\e_1}\al+\fr{\e_2}{\al^2}+\cdots.
$$
Let $\be$ be determined by
$$
1=\fr{\th_{\e_1}}\be+\fr{\th_{\e_2}}{\be^2}+\cdots
$$
Then we have $\al=\al_{\be;\th_0\cdots\th_{q-1}}$. Therefore,
$\al_{\be;\th_1\cdots\th_{q-1}}$ is, and in turn,
$C_{\be;\th_0\cdots\th_{q-1}}$ is continuous with respect to $\be$.
\hskip10pt $\Box$

Now let us have a closer look at $C_{\be;013}$ for $\be\in(3,4]$. By
Theorem \ref{cont} we know that $\dim_H(C_{\be;013})$ is continuous
with respect to $\be$. We have $\dim_H(C_{4;013})=\fr{\log3}{\log4}$
and
$$
\lim_{\be\downarrow3}\dim_H(C_{\be;013})=\fr{\log2}{\log3}=\dim_H(C_{3,01}).
$$
As stated in the above, $\dim_H(C_{\be;013})$ is decreasing for
$\be\in\left(\fr{3+\sqrt{17}}2,\fr{3+\sqrt{21}}2\right)$. In fact,
if $\be_l$ and $\be_r$ be determined by
\begin{equation}
1=\fr3{\be_l}+\sum_{i=2}^{m}\fr{a_i}{\be_l^i}+\fr2{\be_l^{m+1}}
\label{beta11}
\end{equation}
 and
\begin{equation}
1=\fr3{\be_r}+\sum_{i=2}^{m}\fr{a_i}{\be_r^i}+\fr3{\be_r^{m+1}}
\label{beta22}
\end{equation}
where $a_i\in\{0,1,3\}$ with
$$
(a_i,a_{i+1},\cdots, a_m,3, 0, \cdots)<(3,a_1,\cdots,a_m,
3,0,\cdots),
$$
then for any $\be\in[\be_l,\be_r]$ we have the same value for
$\al_{\be;013}$ which is given by
\begin{equation}
1=\fr2{\al^k}+\sum_{i=2}^{m}\fr{b_i}{\al^i}+\fr2{\al^{m+1}}
\label{beta22}
\end{equation}
where $b_i=a_i$ if $a_i\in\{0,1\}$ or $b_i=2$ if $a_i=3$.

Let $V=[3,4]\backslash\cup(\be_l,\be_r)$ where the union is for all
possible $\be_l$, $\be_r$ defined in the above.

\begin{theorem} $\be\in V$ if and only if the $\be$-expansion of 1
is
\begin{equation}
1=\fr{a_1}{\be}+\fr{a_2}{\be^2}+\cdots\label{equ11}
\end{equation}
or
\begin{equation}
1=\fr{a_1}{\be}+\cdots+\fr{a_k}{\be^k}+\fr{2}{\be^{k+1}}\label{equ22}
\end{equation}
where $a_i\in\{0,1,3\}$. For $\be\in(3,4)\backslash V$ we have
$$
\fr{d\dim_\be(C_\be;013)}{d\be}<0.
$$
If $\be\in V$ but $\be\ne\be_l$ or $\be_r$,
$$
\fr{d\dim_H(C_{\be;013})}{d\be}=\infty.
$$
If $\be=\be_l$ or $\be_r$, then the left sided or right sided
derivative of $\dim_H(C_{\be;013})$ with respect to $\be$ is
$\infty$. \label{C013}
\end{theorem}

\noindent{\it Proof.} Given $\be\in(3,4)$, if the $\be$-expansion of
1 contains digit 2, then $\be\in[\be_l,\be_r)$ for some pair of
$\be_l$, $\be_r$ defined in the above. Hence we have either
$\be\notin V$ or $\be=\be_l$. Therefore, $\be\in V$ implies
(\ref{equ11}) or (\ref{equ22}) and vise versa.

Since for $\be\in[\be_1,\be_2]$ we have the same value of
$\al_{\be;013}$, we have $\fr{d\dim_\be(C_\be;013)}{d\be}<0$ for
$\be\in(3,4)\backslash V$.

Taking $\be\in V$, we assume that $\be$ is non-simple and show that
$\lim_{\be'\downarrow\be}\fr{\al-\al'}{\be-\be'}=\infty$. For other
cases the discussion is similar and simpler. Let
$$
1=\fr{a_1}{\be}+\fr{a_2}\be+\cdots
$$
where $a_i\in\{0,1,3\}$ and $a_i\ne0$ for infinitely many $i$. Given
$k\ge1$, choose $i_k> k$ such that $a_{i_k}\ne3$. Let $\be_k$ be
given by
$$
1=\sum_{i=1}^{i_k-1}\fr{a_i}{\be_k^i}+\fr1{\be_k^{i_k-1}}
\left(\fr{a_1}{\be_k}+\fr{a_2}{\be_k^2}+\cdots\right).
$$
Use $\al$ and $\al_k$ to denote $\al_{\be;013}$ and
$\al_{\be_k;013}$. Then
$$
1=\fr{b_1}{\al}+\fr{b_2}{\al^2}+\cdots
$$
and
$$
1=\sum_{i=1}^{i_k-1}\fr{b_i}{\al_k^i}+\fr1{\al^{k_i-1}}
\left(\fr{b_1}{\al_k}+\fr{b_2}{\al_k^2}+\cdots\right).
$$
where $b_i=a_i$ if $a_i=0$ or 1 and $b_i=2$ if $a_i=3$. We estimate
$\be_k-\be$ and $\al_k-\al$.
$$
\begin{aligned}
&\be_k-\be\\
 =&\sum_{i=2}^{i_k-1}\left(\fr{a_i}{\be_k^i}-\fr{a_i}{\be^i}\right)
+\left(\fr{a_1}{\be_k^{i_k}}-\fr{a_{i_k}}{\be^{i_k}}\right)
+\left(\fr{a_{2}}{\be_k^{i_k+1}}-\fr{a_{i_k+1}}{\be^{i_k+1}}\right)+\cdots\\
\le&\fr{a_1-a_{i_k}}{\be_k^{i_k}}+\cdots<\fr c{\be^{i_k}}
\end{aligned}
$$
for some constant $c$.
$$
\begin{aligned}
&\al_k^{i_k}-\al^{i_k}\\
=&\sum_{i=1}^{i_k-1}b_i(\al_k^{i_k-i}-\al^{i_k-i})+(b_1-b_{i_k})+
T_{\al_k} y-T_\al^{i_k}1
\end{aligned}
$$
where $y=\fr{b_1}{\be_k}+\fr{b_2}{\be_k^2}+\cdots$ and
$T_\al^{i_k}1=\fr{b_{i_k+1}}\be+\fr{b_{i_k+2}}{\be^2}+\cdots$.
Noting that $b_1=a_1-1=2$ and $b_{i_k}\le 1$,
$$
\al_k-\al>\fr{1+T_{\al_k}y-T_\al^{n_k}1}{\al_k^{i_k-1}+\al_k^{i_k-1}\al+\cdots+\al^{i_k-1}}
>\fr{T_{\al_k}y}{i_k\al_k^{i_k}}.
$$
Since $\al_k\to\al$ when $k\to\infty$, we have
$$
T_{\al_k}y=\fr{b_2}{\be_k}+\fr{b_3}{\be_k^2}+\cdots\to
\fr{b_2}{\be}+\fr{b_3}{\be^2}=T_\al1>0.
$$
Hence
$$
\fr{\al_k-\al}{\be_k-\be}>\fr{cT_{\al_k}1}{i_k}\left(\fr{\al_k}{\be_k}\right)^{i_k}\to
\infty.
$$
This shows that $
\limsup_{\be'\downarrow\be}\fr{\al_k-\al}{\be_k-\be}=\infty $, which
implies that $\lim_{\be'\downarrow\be}\fr{\al_k-\al}{\be_k-\be}$
since $\al_{\be;013}$ is a monotone function of $\be$.

Now we have shown that for $\be\in V$ if it is non-simple then the
right sided derivative of $\al_{\be;013}$ with respect to $\be$ is
$\infty$. By this we see that the right sided derivative of
$\dim_H(C_{\be;013})$ with respect to $\be$ is $\infty$. By similar
discussions we can obtain that
$$
\fr{d\dim_H(C_{\be;013})}{d\be}=\infty,\ \ \ \be\in V, \text{ but
}\be\ne\be_l\text{ or }\be_r
$$
and
$$
\fr{d\dim_H(C_{\be;013})}{d\be^-}=\infty
$$
if $\be=\be_l$,
$$
\fr{d\dim_H(C_{\be;013})}{d\be^+}=\infty
$$
if $\be=\be_r$,
 where we use ``$-$" and ``$+$" to denote the left and
right handed
 derivatives. \hskip10pt $\Box$

\noindent{\bf Remark.} Although Theorem \ref{C013} is stated and
proved for $C_{\be;013}$, we have similar result for
$C_{\be;\th_0\cdots\th_{q-1}}$ in general. We omit the details of
the statement and proof in general case.

%
%

\section{View $C_{\be;02}$ through Different Eyes.}

In this section, we assume that $2<\be\le3$. Let $F_\be=\{[0,1];
f_0: [0,1]\mapsto [0,1], f_2:[0,\be-2]\mapsto[0,1]\}$ be the local
IFS defined in Section 1,  where $f_0(x)=\fr x\be$ and
$f_2(x)=\fr{x+2}\be$. We show that $C_{\be;02}$ is the unique
invariant set of $T_\be$, for all $\be$ except a countable set.

\begin{theorem} For all $\be\in(2,3]$, $\C$ is an invariant set of $T_\be$.
Let $Q$ be the set of $\be\in(2,3]$ such that the $\be$-expansion of
1 is
$$
1=\fr{\e_1}{\be}+\cdots+\fr{\e_k}{\be^k}
$$
where $\e_i\in\{0,2\}$. For $\be\in(2,3]\backslash Q$, $\C$ is the
only invariant set of $T_\be$ \label{theo4}
\end{theorem}

\noindent{\it Proof.} First we prove that $C_{\be;02}$ is an
invariant set of $F_\be$ by showing that
\begin{equation}
C_{\be;02}=f_0(C_{\be;02})\cup f_2(\C\cap[0,\be-2]).\label{eqqq}
\end{equation}
The compactness of $\C$ ensures that $f_0(\C)\cup
f_2(\C\cap[0,\be-2])$ is compact. For any $x=(0.a_1a_2\cdots)_\be$
with $a_i\in\{0,2\}$ for all $i$, we have $T_\be
x=(0.a_2a_3\cdots)_\be\in \C$. Then $x=f_0(T_\be x)\in f_0(\C)$ if
$a_1=0$, or $T_\be x\in [0,\be-2)$ and $x=f_2(T_\be x)\in f_2(\C\cap
[0,\be-2])$ if $a_1=2$. Hence $f_0(\C)\cup f_2(\C\cap[0,\be-2])$
contains all $x\in[0,1)$ whose $\be$-expansion has only 0 or 2. By
the definition of $\C$ and the compactness of $f_0(\C)\cup
f_2(\C\cap[0,\be-2])$ we get $\C\subset f_0(\C)\cup
f_2(\C\cap[0,\be-2])$. On the other hand, it is easy to see that
$f_0(\C)\subset\C$ and $f_2(\C\cap[0,\be-2])\subset\C$. Therefore,
we have (\ref{eqqq}).

Next show that if $\be\in (2,3]\backslash Q$ then $\C$ is the only
invariant set of $F_\be$. Let $A$ be an invariant set of $T_\be$.
Then $A$ is compact and
\begin{equation}
A=f_0(A)\cup f_2(A\cap[0,\be-2]).\label{poi}
\end{equation}
It is easy to see that $0\in A$. By this and (\ref{poi}) we get that
$A$ contains all $x$ whose $\be$-expansion is finite with all
entries equal to 0 or 2. Then we get $\C\subset A$.

Assume that there exists $x_0\in A\backslash\C$. We claim that we
have $1\in A$. If $x_0\ne1$ then  $x_0=(0.a_1a_2\cdots)_\be$ with
$a_k=1$ for some $k$. Let $k$ be the smallest integer with $a_k=1$.
By (\ref{poi}) we have either $x_0\in f_0(A)$ or $x_0\in
f_2(A\cap[0,\be-2])$. If $x_0\in f_0(A)$, then $y_0=\be x_0\in A$.
This indicates that $a_1=0$ or $x_0=\fr1\be$. If $x_0\in
f_2(A\cap[0,\be-2])$, then $y_0=\be x_0-2\in A$. Hence we always
have $T_\be x_0\in A$. Continue this discussion, we get $T_\be^{k-1}
x_0=(0.{a_k}a_{k+1}\cdots)_\be\in A$. Because $a_k=1$ we have
$T_\be^{k-1}x_0<\fr2\be$. Then $T_\be^{k-1}x_0\notin
f_2(A\cap[0,\be-2]$. If $T_\be^{k-1}x_0>\fr1\be$, then it is not
contained in $f_0(A)$. Hence we must have $T_\be^{k-1}x=\fr1\be$. By
$\fr1\be\in A$ we get $1\in A$.

 If $1\in A$ then, by the above discussion, $T^k1\in A$ for any $k\ge0$.
 By (\ref{poi}) we see that either $\be$-expansion of 1 does not contain 1 or it
 has an only 1 as the last non-zero term. However, if $\be$ is
 non-simple and the $\be$-expansion of 1 does not contain 1, or $\be$ is
 simple with the only digit 1 as the last non-zero term, then we have
 $1\in \C$ which implies $x_0\in\C$. Hence only when $\be$ is simple
and the $\be$-expansion of does not contain 1, that is, $\be\in Q$,
we have $1\in A\backslash \C$. Therefore, $\C$ is the only invariant
set of $T_\be$ for $\be\in(2,3]\backslash Q$. \hskip10pt $\Box$

Construct a compact set $B_{\be;02}$ by a interval removal process
similar to that for the classical Cantor set. Firstly we remove the
interval $(\frac1\be,\fr2{\be})$ from the unit interval [0,1]. Next
we remove $(\frac1{\be^2},\frac2{\be^2})$ from the remaining
interval $[0,\fr1\be]$ and remove
$(\frac2\be+\frac1{\be^2},\frac2\be+\fr2{\be^2})\cap[\frac2\be,1]$.
from $[\fr2\be,1]$. In general, at step $n$ if a remaining interval
$[u_n,v_n]$ has $v_n-u_n>\fr1{\be^n}$, we remove
$(u_n+\fr1{\be^n},u_n+\fr2{\be^n})\cap[u_n,v_n]$ from it.

\begin{theorem} $B_{\be;02}$ is an invariant set of $T_\be$ for all $\be\in(2,3]$.
For $\be\in Q$, we have $B_{\be;02}\ne\C$ and
$B_{\be;02}\backslash\C$ is a countable set of isolated
points.\label{theo6}
\end{theorem}

\noindent{\it Proof}. First we show that $B_{\be;02}=\C$ for
$\be\in(2,3]\backslash Q$. It is clear that $B_{\be;02}\supset\C$.
By the construction process, if the $\be$-expansion of a number has
a 1 and at least one non-zero term after the digit 1 then it will be
removed at some stage. Hence $B_{\be;02}$ consists of points in $\C$
and possibly the numbers whose $\be$-expansion is finite with an
only digit 1 as the last non-zero term. As discussed above, when
$\be\in(2,3]\backslash Q$ these numbers are contained in $\C$. Hence
we have $B_{\be,02}=\C$ for $\be\in(2,3]\backslash Q$.

Next we show that for $\be\in Q$ we have $B_{\be;02}\ne\C$ and
$B_{\be;02}$ is the only invariant set of $T_\be$ other than $\C$.
Assume that $1=(0.\e_1\e_2\cdots\e_{k-1}2)_\be$. By the construction
of $B_{\be;02}$, the interval
$\left(\fr{\e_1}{\be}+\cdots+\fr{\e_{k-1}}{\be^{k-1}}+\fr1{\be^k},
\fr{\e_1}{\be}+\cdots+\fr{\e_{k-1}}{\be^{k-1}}+\fr2{\be^k}\right)$
is removed and 1 remains in $B_{\be;02}$. Hence $B_{\be;02}\ne\C$.
It is easy to see that in this case $B_{\be;02}$ contains all
$x=(0.a-1\cdots a_{k-1}1)_\be$ where $a_i\in\{0,2\}$. Therefore
$B_{\be;02}$ is just the invariant set $A$ in the proof of
(\ref{theo4}) when $A\ne\C$. Hence $B_{\be;02}$ is an invariant set
of $T_\be$ and $T_\be$ does not have invariant set other than
$B_{\be;02}$ and $\C$. It is obvious that $B_{\be;02}\backslash\C$
consists of countablly many isolated points.\hskip10pt $\Box$

We can also define a normal IFS $([0,1]; f_0, f_2)$  by extending
$f_2$ to the whole interval with
$$
f_2(x)=\begin{cases} \fr{x+2}\be &\ x\in[0,\be-2],\\
                              1     &\ x\in(\be-2,1].
\end{cases}
$$
Then $B_{\be;02}$ is the attractor of this IFS. We omit the details.



\section{Further Research.}

1. In Theorem 1 we obtained the Hausdorff dimension for
$\be$-attractors when all the maps have the same contract ratio
under a separation condition\break $f_i(E_\be)\cap f_j(E_\be)\cap
E_\be=\emptyset$, whenever $i\ne j$. We believe that the separation
condition can be replaced by an open set condition. Comparing with
the $\be$-expansion of real numbers in $[0,1)$, we may define a
$\be$-open set condition as follows.

{\bf $\be$-Open Set Condition.} Let $F=(X; f_0,f_1,\cdots, f_{n-1})$
be an IFS where $X$ is a compact subset of $R^m$ and $f_i$ is a
Lipschiz map with Lipschiz constant $r_i<1$ for $0\le i\le n-1$. Let
$n-1<\be\le n$ and $E_\be$ is the $\be$-attractor. We say $F$
satisfies a $\be$-open set condition if there exists an open set $O$
with $E_\be\subset \bar{O}$ such that $f_i(O)\subset O$ and
$f_i(O)\cap f_j(O)=\emptyset$ for $0\le i<j\le n-1$.

Note that in the above definition we do not require $f_{n-1}(O)\subset O$
but assume $E_\be\subset \bar{O}$ 
which is automatically true in the original version of open set
condition. It is easy to see that $\CC$ satisfies a $\be$-open set
condition if we choose $O=(0,1)$. When $\be=n$ is an integer the
original version of open set condition implies the $\be$-open set
condition for the same open set $O$. When all $f_i$'s are
similarities, a $\be$-open set condition with $\be$ being an integer
also implies the original version of open set condition possibly
with a different open set $O'$.

\noindent{\bf Conjecture.} {\it Let $(X; f_0,f_1,\cdots,f_{n-1})$ be
an IFS such that
 $$
 d(f_i(x), f_i(y))=rd(x,y)
 $$
 for all $x, y\in X$ and $0\le i\le n-1$, where $0<r<1$.
Let $1<\be<n$ and $E_\be$ be the $\be$-attractor. Under the
$\be$-open set condition, we have
$$
\dim(E_\be)=\frac{\log\be}{-\log r}.
$$
The Hausdorff measure of $E_\be$ in its dimension is positive and
finite. }


2. Given an IFS $(X; f_0,\cdots, f_{n-1})$ with attractor $E$ and a
probability vector $(p_0,\cdots, p_{n-1})$ with $p_i>0$, for any
Borel measure $\mu$ on $X$ we define
$$
F(\mu)(B)=\sum_{i=0}^{n-1}p_i\mu (f_i^{-1}(B))
$$
for any Borel subset $B$. Then there is a unique measure $\nu$ such
that $F(\nu)=\nu$. $\nu$ has $E$ as its support and for any $\mu$
the sequence $\{F^{(k)}(\mu)\}$ weakly converges to $\nu$.

 It is natural to study the invariant measure for $\be$-attractors.

 \noindent{\bf Acknowledgement.} This work was supported
by an Australian Research Council discovery grant. The author is
grateful to Professor John Hutchinson for many discussions with him.
The author thanks Professor Michael Barnsley for some helpful
comments.


\begin{thebibliography}{99}
\bibitem{Aki1} S. Akiyama and N. Gjini, {\it On the connectedness of self-affine
attractors}, Arch. Math. (Basel) {\bf 82} (2004), 153–-163.

\bibitem{Aki2} S. Akiyama and T. Sadahiro, Taizo, {\it A self-similar tiling
generated by the minimal Pisot number}, Proceedings of the 13th
Czech and Slovak International Conference on Number Theory
(Ostravice, 1997). Acta Math. Inform. Univ. Ostraviensis {\bf 6}
(1998), 9–-26.


\bibitem{bandt} C. Bandt, {\it Self-similar Sets 1. Markov Shifts and Mixed Self-similar Sets,}
Math. Nachr. {\bf 142} (1989), 107-123.

\bibitem{bar} Luis M. Barreira, {\it Cantor Sets with Complicated
Geometry and Modeled by General Symbolic Dynamical Systems}, Rand.
Comp. Dynam. {\bf 3} (1995), 213-139.

\bibitem{BM} A. Bertrand-Mathis, {\it Questions Diverses Relative
aux Syst\`ems cod\'es: Applications au $\theta$-shift}, Preprint.

\bibitem{mbarnsley} Michael F. Barnsley and Lyman P. Hurd, Fractal
Image Compression, AK Peters Ltd., Wellesley, Massachusetts, 1993.
\bibitem{barnsley} Michael F. Barnsley, {\it Theory and application of
fractal tops}, 3—20, Fractals in Engineering: New Trends in Theory
and Applications. Lévy-Véhel J.; Lutton, E. (eds.) Springer-Verlag,
London Limited, 2005.

\bibitem{Brown-Yin2} Gavin Brown and Qinghe Yin, {\it $\beta$-transformation,
natural extension and invariant measure},  Ergodic Theory Dynam.
Systems {\bf 20}  (2000), 1271--1285.

\bibitem{Falconer} K.J. Falconer, Fractal Geometry, John Wiley, 1990.
\bibitem{Falconer97}K.J. Falconer, Techniques in Fractal Geometry, John Wiley, 1997.

\bibitem {Fan-Zhu} Aihua Fan and Hao Zhu, {\it Level sets of
$\beta$-expansions},  Comptes Rendus Mathematique {\bf 339} (2004),
709--712.

\bibitem{Hutchinson} John E. Hutchinson, {\it Fractals
and self-similarity}, Indiana Univ. Math. J. {\bf 30} (1981),
713--747.

\bibitem{Ito} S. Ito and M. Kimura, {\it On Rauzy fractal}, Janpan
J. Indust.Appl. Math. {\bf 8} (1991), 461--486.

\bibitem{Ito2} S. Ito and Y. Sano, {\it On periodic $\beta$-expansions
of Pisot numbers and Rauzy fractals}, Osaka J. Math. {\bf 38}
(2001),  349–-368.

\bibitem{keaneetal}M. Keane, M. Smorodinsky and B. Solomyk, {On the
Morphology of $\gamma$-expansions  with Deleted Digits}, Thans. Amer. Math. Soc.
{\bf 347} (1995), 955-966.

\bibitem{lind} D. Lind and B. Marcus, An Introduction to Symbolic
Dynamical Systems and Coding, Cambridge University Press, 1995.

\bibitem{Luz} D. Luzeaux, {\it From $\beta$-expansions to chaos and
fractals}, Complexity Internat. {\bf 1} (1994).

\bibitem{mauldin}R.D. Mauldin and S.C. Williams, {\it Hausdorff
dimension in graph directed constructions}, Trans. Amer. Math. Soc.
{\bf 309} (1988), 811-829.

\bibitem{Mori}M. Mori, {\it Piecewise linear transformations of the
 unit interval and Cantor sets}, in N. Pythyeas Fogg, Substitutions
 in Dynamics, Arithemetics and Combinatorics (Lecture Notes in
 Mathematics, 1794), Springer-Verlag Berlin 2002, 343--362.

\bibitem{Parry}W. Parry, {\it On the $\be$-expansions of real numbers},
 Acta Math. Hungar. {\bf 11} (1960), 401--416.

\bibitem{pesin} Y. Pesin and H. Weiss, {\it On the dimension
of deterministic and random Cantor-like sets, symbolic dynamics, and
the Eckmann-Ruelle conjecture}, Commun. Math. Phys. {\bf 182}
(1996), 105-153.

\bibitem{pollicott}M. Pollicott and K. Simon, {\it The Hausdorff Dimension
of $\lambda$-expansions with deleted digits}, Trans. Amer. Math. Soc.
{\bf 347} (1995), 967-943

\bibitem{Renyi}A. R\'enyi, {\it Representations for real numbers and their
ergodic properties}, Acta Math. Hungar. {\bf 8} (1957), 477--493.

\bibitem{ru} D. Ruelle, Thermodynamics Formalism, Reading,
Addison-Wesley, 1978.

\bibitem{sch} J. Schmeling, {\it Most $\beta$-shifts Have Bad
Ergodic Properties}, Preprint, IAAS (1993).

\bibitem{solomyak}B. Solomyak, {\it Measure and Dimension for some Fractal Families},
Math. Proc. Camb. Phil. Soc. {\bf 124} (1998), 531-546.

\bibitem{yinq}Qinghe Yin, {\it On Hausdorff dimension for attractors
of iterated function systems}, J. Austal. Math. Soc. (Series A) {\bf
55} (1993), 216-231.
\end{thebibliography}
\end{document}